\documentclass[12pt]{article}

\usepackage{amsmath, epsfig, cite}
\usepackage{amssymb}
\usepackage{amsfonts}
\usepackage{latexsym}
\usepackage{float}

\newtheorem{thm}{Theorem}[section]

\newtheorem{conj}[thm]{Conjecture}
\newtheorem{lem}[thm]{Lemma}

\numberwithin{equation}{section}

\newcommand{\qed}{{\hfill$\square$}\medskip}

\newcommand*{\pFq}[5]{{}_{#1}F_{#2}\left[ \begin{matrix} #3\\[5pt] #4\end{matrix};#5\right]}

\begin{document}

\begin{center}
{\Large\bf On van Hamme's (A.2) and (H.2) supercongruences}
\end{center}

\vskip 2mm \centerline{Ji-Cai Liu}
\begin{center}
{\footnotesize Department of Mathematics, Wenzhou University, Wenzhou 325035, PR China\\
{\tt jcliu2016@gmail.com} }
\end{center}


\vskip 0.7cm \noindent{\bf Abstract.}
In 1997, van Hamme conjectured 13 Ramanujan-type supercongruences labeled (A.2)--(M.2).
Using some combinatorial identities discovered by {\tt Sigma}, we extend (A.2) and (H.2) to supercongruences modulo $p^4$ for primes $p\equiv 3\pmod{4}$, which appear to be new.

\vskip 3mm \noindent {\it Keywords}:
Supercongruences; Hypergeometric series; $p$-Adic Gamma functions
\vskip 2mm
\noindent{\it MR Subject Classifications}: 05A19, 11A07, 33C20
\section{Introduction}
In 1913, Ramanujan announced the following identity in his famous letter:
\begin{align}
\sum_{k=0}^{\infty}(-1)^k(4k+1)\left(\frac{\left(\frac{1}{2}\right)_k}{k!}\right)^5
=\frac{2}{\Gamma\left(\frac{3}{4}\right)^4},\label{new-1}
\end{align}
where $(a)_0=1$ and $(a)_n=a(a+1)\cdots (a+n-1)$ for $n\ge 1$. The formula \eqref{new-1} was later proved by Hardy (1924) and Watson (1931).

Motivated by some of such formulas, van Hamme \cite{vanHamme-1997} conjectured  13 supercongruences in 1997, which relate the partial sums of certain hypergeometric series to the values of the $p$-adic Gamma functions.
These 13 conjectural supercongruences labeled (A.2)--(M.2) are $p$-adic analogues of
their corresponding infinite series identities. For instance, van Hamme conjectured that the identity \eqref{new-1} has the following interesting $p$-adic analogue, which was first proved by McCarthy and Osburn \cite{mo-2008} by use of Gaussian hypergeometric series.
\begin{thm} (See \cite[(A.2)]{vanHamme-1997}.)
For any odd prime $p$, we have
\begin{align*}
\sum_{k=0}^{\frac{p-1}{2}}(-1)^k(4k+1)\left(\frac{\left(\frac{1}{2}\right)_k}{k!}\right)^5\equiv
\begin{cases}
-p\Gamma_p\left(\frac{1}{4}\right)^4\pmod{p^3}\quad &\text{if $p\equiv 1\pmod{4}$,}\\[10pt]
0\pmod{p^3}\quad &\text{if $p\equiv 3\pmod{4}$,}
\end{cases}
\end{align*}
where $\Gamma_p(\cdot)$ denotes the $p$-adic Gamma function recalled in the next section.
\end{thm}

Later, Swisher \cite{swisher-2015} showed that the supercongruence (A.2) also holds modulo $p^5$ for primes $p\equiv 1\pmod{4}$. It is natural to ask whether the case $p\equiv 3\pmod{4}$ possesses a modulo $p^5$ extension in terms of the $p$-adic Gamma functions. Interestingly, numerical calculation suggests the following result which partially motivates this paper.

\begin{conj}
For primes $p\ge 5$ with $p\equiv 3\pmod{4}$, we have
\begin{align}
\sum_{k=0}^{\frac{p-1}{2}}(-1)^k(4k+1)\left(\frac{\left(\frac{1}{2}\right)_k}{k!}\right)^5\equiv
-\frac{p^3}{16}\Gamma_p\left(\frac{1}{4}\right)^4\pmod{p^5}.\label{ss-1}
\end{align}
\end{conj}

It is very difficult for us to prove \eqref{ss-1} thoroughly, but we can show that \eqref{ss-1} holds modulo $p^4$.
\begin{thm}\label{t1}
Let $p\ge 5$ be a prime. For $p\equiv 3\pmod{4}$, we have
\begin{align}
\sum_{k=0}^{\frac{p-1}{2}}(-1)^k(4k+1)\left(\frac{\left(\frac{1}{2}\right)_k}{k!}\right)^5\equiv
-\frac{p^3}{16}\Gamma_p\left(\frac{1}{4}\right)^4\pmod{p^4}.\label{aa3}
\end{align}
\end{thm}

The other motivating example of this paper is the following supercongruence which was
both conjectured and confirmed by van Hamme \cite{vanHamme-1997}.
\begin{thm} (See \cite[(H.2)]{vanHamme-1997}.)
For any odd prime $p$, we have
\begin{align*}
\sum_{k=0}^{\frac{p-1}{2}}\left(\frac{\left(\frac{1}{2}\right)_k}{k!}\right)^3\equiv
\begin{cases}
-\Gamma_p\left(\frac{1}{4}\right)^4\pmod{p^2}\quad &\text{if $p\equiv 1\pmod{4}$,}\\[10pt]
0\pmod{p^2}\quad &\text{if $p\equiv 3\pmod{4}$.}
\end{cases}
\end{align*}
\end{thm}

Recently, Long and Ramakrishna \cite[Theorem 3]{lr-2016} obtained a modulo $p^3$ extension of (H.2) as follows:
\begin{align}
\sum_{k=0}^{\frac{p-1}{2}}\left(\frac{\left(\frac{1}{2}\right)_k}{k!}\right)^3\equiv
\begin{cases}
-\Gamma_p\left(\frac{1}{4}\right)^4 \pmod{p^3}\quad &\text{if $p\equiv 1\pmod{4}$,}\\[10pt]
-\frac{p^2}{16}\Gamma_p\left(\frac{1}{4}\right)^4\pmod{p^3}\quad &\text{if $p\equiv 3\pmod{4}$.}
\end{cases}\label{aa2}
\end{align}

The second aim of this paper is to establish a modulo $p^4$ extension of \eqref{aa2} for primes
 $p\equiv 3\pmod{4}$.
\begin{thm}\label{t2}
Let $p\ge 5$ be a prime. For $p\equiv 3\pmod{4}$, we have
\begin{align}
\sum_{k=0}^{\frac{p-1}{2}}\left(\frac{\left(\frac{1}{2}\right)_k}{k!}\right)^3
\equiv -\frac{p^2}{4}\Gamma_p\left(\frac{1}{4}\right)^4J_{(p-3)/4}^2 \pmod{p^4},\label{aa4}
\end{align}
where $J_n=\sum_{k=0}^n{2k\choose k}^2/16^k$.
\end{thm}

{\noindent \bf Remark.}
We shall show that $J_{(p-3)/4}\equiv -\frac{1}{2}\pmod{p}$, which implies \eqref{aa4}
is indeed a modulo $p^4$ extension of the second case of \eqref{aa2}.
By the Chu-Vandermonde identity, we have
\begin{align*}
\sum_{k=0}^n{2n+1\choose k}^2=\frac{1}{2}\sum_{k=0}^{2n+1}{2n+1\choose k}^2=\frac{1}{2}{4n+2\choose 2n+1}.
\end{align*}
Letting $n=\frac{p-3}{4}$ with $p\equiv 3\pmod{4}$ in the above and noting that
\begin{align*}
{\frac{p-1}{2}\choose k}\equiv \frac{{2k\choose k}}{(-4)^k}\pmod{p}
\quad \text{and}\quad {p-1\choose \frac{p-1}{2}}\equiv -1\pmod{p},
\end{align*}
we immediately obtain $J_{(p-3)/4}\equiv -\frac{1}{2}\pmod{p}$.

It is known that Gaussian hypergeometric series (see, for example, \cite{Long-2011,lr-2016,mo-2008,mortenson-2008,swisher-2015}) and the W-Z method (see, for example, \cite{guojw-itsf-2017,oz-2016,sunzw-2012,zudilin-2009}) commonly apply to the Ramanujan-type supercongruences. We refer to \cite{oz-2016,swisher-2015} for more recent developments on van Hamme's supercongruences. In this paper, we provide a different approach which is based on some combinatorial identities involving harmonic numbers. All of these identities are automatically discovered and proved by the software package {\tt Sigma} developed by Schneider \cite{schneider-1999}.

The rest of this paper is organized as follows. The preliminary section, Section 2, is devoted to some properties of the $p$-adic Gamma function. We prove Theorems \ref{t1} and \ref{t2} in Sections 3 and 4, respectively.

\section{Preliminary results}
We first recall the definition and some basic properties of the $p$-adic Gamma function. For more details, we refer to \cite[Section 11.6]{cohen-2007}.
Let $p$ be an odd prime and $\mathbb{Z}_p$ denote the set of all $p$-adic integers. For
$x\in \mathbb{Z}_p$, the $p$-adic Gamma function is defined as
\begin{align*}
\Gamma_p(x)=\lim_{m\to x}(-1)^m\prod_{\substack{0< k < m\\
(k,p)=1}}k,
\end{align*}
where the limit is for $m$ tending to $x$ $p$-adically in $\mathbb{Z}_{\ge 0}$.

\begin{lem} (See  \cite[Section 11.6]{cohen-2007}.)
For any odd prime $p$ and $x\in \mathbb{Z}_p$, we have
\begin{align}
&\Gamma_p(1)=-1,\label{bb1}\\
&\Gamma_p(x)\Gamma_p(1-x)=(-1)^{s_p(x)},\label{bb2}\\
&\frac{\Gamma_p(x+1)}{\Gamma_p(x)}=
\begin{cases}
-x\quad&\text{if $|x|_p=1$,}\\
-1\quad &\text{if $|x|_p<1$, }
\end{cases}\label{bb3}
\end{align}
where $s_p(x)\in \{1,2,\cdots,p\}$ with $s_p(x)\equiv x\pmod{p}$ and
$|\cdot|_p$ denotes the $p$-adic norm.
\end{lem}

\begin{lem} (See \cite[Lemma 17, (4)]{lr-2016}.)
Let $p$ be an odd prime. If $a\in \mathbb{Z}_p, n\in \mathbb{N}$
such that none of $a,a+1,\cdots,a+n-1$ in $p\mathbb{Z}_p$, then
\begin{align}
(a)_n=(-1)^n\frac{\Gamma_p(a+n)}{\Gamma_p(a)}.\label{bb4}
\end{align}
\end{lem}

The following lemma is a special case of a theorem due to Long and Ramakrishna \cite{lr-2016}.
\begin{lem} (See \cite[Theorem 14]{lr-2016}.)
For any prime $p\ge 5$ and $a,b\in\mathbb{Z}_p$, we have
\begin{align}
\Gamma_p(a+bp)\equiv \Gamma_p(a)\left(1+G_1(a)bp\right) \pmod{p^2},\label{bb5}
\end{align}
where $G_1(a)=\Gamma_p^{'}(a)/\Gamma_p(a)\in \mathbb{Z}_p$.
\end{lem}

\section{Proof of Theorem \ref{t1}}
In order to prove Theorem \ref{t1}, we establish the following two combinatorial identities.
\begin{lem}
For any odd integer $n$, we have
\begin{align}
\sum_{k=0}^{n}(-1)^k(4k+1)
\frac{\left(\frac{1}{2}\right)_k^3(-n)_k(n+1)_k}{(1)_k^3(n+\frac{3}{2})_k(-n+\frac{1}{2})_k}=0.\label{cc1}
\end{align}
\end{lem}
{\noindent \it Proof.}
We begin with the following identity \cite[(2.1)]{mortenson-2008}:
\begin{align*}
&\pFq{6}{5}{a,&\frac{a}{2}+1,&b,&c,&d,&e}{&\frac{a}{2},&a-b+1,&a-c+1,&a-d+1,&a-e+1}{-1}\\
&\hskip 2cm=\frac{\Gamma(a-d+1)\Gamma(a-e+1)}{\Gamma(a+1)\Gamma(a-d-e+1)}
\pFq{3}{2}{a-b-c+1,&d,&e}{&a-b+1,&a-c+1}{1}.
\end{align*}
Letting $a=d=e=\frac{1}{2},b=-n$ and $c=n+1$ in the above gives
\begin{align}
&\sum_{k=0}^{n}(-1)^k(4k+1)
\frac{\left(\frac{1}{2}\right)_k^3(-n)_k(n+1)_k}{(1)_k^3(n+\frac{3}{2})_k(-n+\frac{1}{2})_k}\notag\\
&\hskip 1cm=
\pFq{3}{2}{\frac{1}{2},&\frac{1}{2},&\frac{1}{2}}{&n+\frac{3}{2},&-n+\frac{1}{2}}{1}{\bigg/}
\Gamma\left(\frac{1}{2}\right)\Gamma\left(\frac{3}{2}\right).\label{cc2}
\end{align}
Recall the following Whipple's identity \cite[page 54]{slater-1966}:
\begin{align}
&\pFq{3}{2}{a,&1-a,&c}{&e,&1+2c-e}{1}\notag\\
&=
\frac{2^{1-2c}\pi\Gamma(e)\Gamma(1+2c-e)}{\Gamma\left(\frac{1}{2}e+\frac{1}{2}a\right)
\Gamma\left(\frac{1}{2}+\frac{1}{2}e-\frac{1}{2}a\right)\Gamma\left(\frac{1}{2}+c-\frac{1}{2}e+\frac{1}{2}a\right)
\Gamma\left(1+c-\frac{1}{2}e-\frac{1}{2}a\right)}.\label{cc3}
\end{align}
Assume that $n$ is odd. Letting $a=c=\frac{1}{2}$ and $e=n+\frac{3}{2}$ in \eqref{cc3} and noting that both $1-\frac{1}{2}(n+1)$ and $\frac{1}{2}(1-n)$ are non-positive integers,
we conclude that
\begin{align}
\pFq{3}{2}{\frac{1}{2},&\frac{1}{2},&\frac{1}{2}}{&n+\frac{3}{2},&-n+\frac{1}{2}}{1}=0.\label{cc4}
\end{align}
Combining \eqref{cc2} and \eqref{cc4}, we complete the proof of \eqref{cc1}. \qed

\begin{lem}
For any odd integer $n$, we have
\begin{align}
&\sum_{k=0}^{n}(-1)^k(4k+1)
\frac{\left(\frac{1}{2}\right)_k^3(-n)_k(n+1)_k}{(1)_k^3(n+\frac{3}{2})_k(-n+\frac{1}{2})_k}
\sum_{j=1}^k\left(\frac{1}{(2j)^2}-\frac{1}{(2j-1)^2}\right)\notag\\
&\hskip 1cm=4^{n-2}(2n+1){\bigg/}n^2{n-1\choose \frac{n-1}{2}}^2.\label{cc5}
\end{align}
\end{lem}
{\noindent \it Proof.}
Actually, we can automatically discover and prove \eqref{cc5} by use of the software package {\tt Sigma} developed by Schneider (see \cite[Section 5]{os-mc-2009} and \cite[Section 3.1]{schneider-1999} for a similar approach to finding and proving identities of this type).

After loading {\tt Sigma} into Mathematica, we insert:\\[7pt]
$\displaystyle \text{\sf In[1]}:=\text{mySum}=\sum_{k=0}^{2n+1}(-1)^k(4k+1)
\frac{\left(\frac{1}{2}\right)_k^3(-2n-1)_k(2n+2)_k}{(1)_k^3(2n+\frac{5}{2})_k(-2n-\frac{1}{2})_k}
\sum_{j=1}^k\left(\frac{1}{(2j)^2}-\frac{1}{(2j-1)^2}\right)$\\[7pt]
We compute the recurrence for the above sum:\\[5pt]
$\text{\sf In[2]}:={\text{rec}=\text{GenerateRecurrence}[\text{mySum},n][[1]]}$\\[5pt]
We find that the above sum $\text{SUM}[n]$ satisfies a recurrence of order $2$. Now we solve this recurrence:\\[5pt]
$\text{\sf In[3]}:=\text{recSol}=\text{SolveRecurrence}[\text{rec},\text{SUM}[n]]$\\[7pt]
{\scriptsize $
\displaystyle\text{\sf Out[3]}=\left\{\left\{0,-\frac{4 (3+4 n)}{(1+2
n)^2}\prod_{\iota _1=1}^n\frac{4 \iota _1^2}{\left(-1+2 \iota _1\right)^2}\right\},
\left\{0,\frac{ 3+4 n}{(1+2 n)^2}\left(\prod_{\iota _1=1}^n\frac{4 \iota _1^2}{\left(-1+2 \iota _1\right)^2}\right)
\left(\sum_{\iota _1=1}^n\frac{4}{\iota _1^2}-\sum_{\iota _1=1}^n\frac{16}{\left(1+2 \iota_1\right){}^2}\right)\right\},\{1,0\}\right\}$}\\[7pt]
Finally, we combine the solutions to represent mySum as follows:\\[5pt]
$\text{\sf In[4]}:=\text{FindLinearCombination}[\text{recSol},\text{mySum},2]$\\[7pt]
$
\displaystyle\text{\sf Out[4]}=\frac{3+4 n}{4 (1+2 n)^2} \prod_{\iota _1=1}^n\frac{4 \iota _1^2}{\left(-1+2 \iota _1\right)^2}
$\\[7pt]
This implies that both sides of \eqref{cc5} are equal.
\qed

{\noindent \it Proof of Theorem \ref{t1}.}
Assume that $p\equiv 3\pmod{4}$. Letting $n=\frac{p-1}{2}$ in \eqref{cc1} yields
\begin{align}
\sum_{k=0}^{\frac{p-1}{2}}(-1)^k(4k+1)
\frac{\left(\frac{1}{2}\right)_k^3(\frac{1-p}{2})_k(\frac{1+p}{2})_k}{(1)_k^3(1+\frac{p}{2})_k(1-\frac{p}{2})_k}
=0.\label{cc6}
\end{align}
We next determine the following product modulo $p^4$:
\begin{align*}
\frac{(\frac{1-p}{2})_k(\frac{1+p}{2})_k}{(1+\frac{p}{2})_k(1-\frac{p}{2})_k}
=\prod_{j=1}^k\frac{(2j-1)^2-p^2}{(2j)^2-p^2}.
\end{align*}
From the following two Taylor expansions:
\begin{align*}
\frac{(2j-1)^2-x^2}{(2j)^2-x^2}=\left(\frac{2j-1}{2j}\right)^2-\frac{4j-1}{(2j)^4}x^2+\mathcal{O}(x^4),
\end{align*}
and
\begin{align*}
\prod_{j=1}^k(a_j+b_jx^2)=\prod_{j=1}^ka_j\cdot\left(1+x^2\sum_{j=1}^k\frac{b_j}{a_j}\right)
+\mathcal{O}(x^4),
\end{align*}
we deduce that for $0\le k\le \frac{p-1}{2}$,
\begin{align}
\frac{(\frac{1-p}{2})_k(\frac{1+p}{2})_k}{(1+\frac{p}{2})_k(1-\frac{p}{2})_k}&\equiv \prod_{j=1}^k\left(\left(\frac{2j-1}{2j}\right)^2-\frac{(4j-1)p^2}{(2j)^4}\right)\notag\\
&\equiv \left(\frac{\left(\frac{1}{2}\right)_k}{(1)_k}\right)^2
\left(1+p^2\sum_{j=1}^k\left(\frac{1}{(2j)^2}-\frac{1}{(2j-1)^2}\right)\right)\pmod{p^4}.\label{cc8}
\end{align}
Substituting \eqref{cc8} into \eqref{cc6} gives
\begin{align}
&\sum_{k=0}^{\frac{p-1}{2}}(-1)^k(4k+1)\left(\frac{\left(\frac{1}{2}\right)_k}{k!}\right)^5\notag\\
&\hskip 1cm\equiv -p^2\sum_{k=0}^{\frac{p-1}{2}}(-1)^k(4k+1)\left(\frac{\left(\frac{1}{2}\right)_k}{k!}\right)^5
\sum_{j=1}^k\left(\frac{1}{(2j)^2}-\frac{1}{(2j-1)^2}\right)\pmod{p^4}.\label{cc9}
\end{align}
Furthermore, letting $n=\frac{p-1}{2}$ in \eqref{cc5} and noting the fact that
\begin{align}
(a+bp)_k(a-bp)_k\equiv (a)_k^2\pmod{p^2},\label{cc10}
\end{align}
we obtain
\begin{align}
&\sum_{k=0}^{\frac{p-1}{2}}(-1)^k(4k+1)\left(\frac{\left(\frac{1}{2}\right)_k}{k!}\right)^5
\sum_{j=1}^k\left(\frac{1}{(2j)^2}-\frac{1}{(2j-1)^2}\right)\notag\\
&\hskip 1cm\equiv p2^{p-5}{\bigg/}\left(\frac{p-1}{2}\right)^2{\frac{p-3}{2}\choose \frac{p-3}{4}}^2\pmod{p^2}.\label{cc11}
\end{align}
It follows from \eqref{cc9} and \eqref{cc11} that
\begin{align}
\sum_{k=0}^{\frac{p-1}{2}}(-1)^k(4k+1)\left(\frac{\left(\frac{1}{2}\right)_k}{k!}\right)^5
\equiv -p^32^{p-5}{\bigg/}\left(\frac{p-1}{2}\right)^2{\frac{p-3}{2}\choose \frac{p-3}{4}}^2\pmod{p^4}.\label{cc12}
\end{align}

By the Fermat's little theorem, we have
\begin{align}
2^{p-5}{\bigg/}\left(\frac{p-1}{2}\right)^2\equiv \frac{1}{4}\pmod{p}.\label{cc13}
\end{align}
On the other hand, using \eqref{bb1}, \eqref{bb4} and \eqref{bb5} we obtain
\begin{align}
{\frac{p-3}{2}\choose \frac{p-3}{4}}^2=2^{p-3}\frac{\left(\frac{1}{2}\right)_{\frac{p-3}{4}}^2}{(1)_{\frac{p-3}{4}}^2}
=2^{p-3}\frac{\Gamma_p\left(\frac{p-1}{4}\right)^2\Gamma_p(1)^2}
{\Gamma_p\left(\frac{1}{2}\right)^2\Gamma_p\left(\frac{p+1}{4}\right)^2}
\equiv \frac{1}{4}\cdot\frac{\Gamma_p\left(-\frac{1}{4}\right)^2}
{\Gamma_p\left(\frac{1}{2}\right)^2\Gamma_p\left(\frac{1}{4}\right)^2}\pmod{p}.\label{cc14}
\end{align}
From \eqref{cc12}--\eqref{cc14}, we deduce that
\begin{align*}
\sum_{k=0}^{\frac{p-1}{2}}(-1)^k(4k+1)\left(\frac{\left(\frac{1}{2}\right)_k}{k!}\right)^5
\equiv -p^3\frac{\Gamma_p\left(\frac{1}{2}\right)^2\Gamma_p\left(\frac{1}{4}\right)^2}{\Gamma_p\left(-\frac{1}{4}\right)^2}
\pmod{p^4}.
\end{align*}
Since $p\equiv 3\pmod{4}$, by \eqref{bb2} we have
\begin{align*}
\Gamma_p\left(\frac{1}{2}\right)^2=(-1)^{\frac{p+1}{2}}=1
\quad\text{and}\quad
\Gamma_p\left(-\frac{1}{4}\right)^2\Gamma_p\left(\frac{5}{4}\right)^2=1.
\end{align*}
It follows that
\begin{align}
\sum_{k=0}^{\frac{p-1}{2}}(-1)^k(4k+1)\left(\frac{\left(\frac{1}{2}\right)_k}{k!}\right)^5
\equiv -p^3\Gamma_p\left(\frac{1}{4}\right)^2\Gamma_p\left(\frac{5}{4}\right)^2
\pmod{p^4}.\label{cc15}
\end{align}
Finally, by \eqref{bb3} we obtain
\begin{align*}
\Gamma_p\left(\frac{5}{4}\right)=-\frac{1}{4}\Gamma_p\left(\frac{1}{4}\right).
\end{align*}
Substituting the above into \eqref{cc15}, we complete the proof of \eqref{aa3}.
\qed

\section{Proof of Theorem \ref{t2}}
\begin{lem}
For any odd integer $n$, we have
\begin{align}
\sum_{k=0}^n\frac{\left(\frac{1}{2}\right)_k(-n)_k(n+1)_k}{(1)_k^3}=0.\label{dd1}
\end{align}
\end{lem}
{\noindent \it Proof.}
Assume that $n$ is odd. Letting $a=-n,e=1$ and $c=\frac{1}{2}$ in \eqref{cc3} and noting that
$\frac{1-n}{2}$ is a non-positive integer, we immediately arrive at \eqref{dd1}.
\qed

\begin{lem}
For any odd integer $n$, we have
\begin{align}
\sum_{k=0}^n\frac{\left(\frac{1}{2}\right)_k(-n)_k(n+1)_k}{(1)_k^3}\sum_{j=1}^k\frac{1}{(2j-1)^2}
=-4^{n-1}\left(\sum_{k=0}^{\frac{n-1}{2}}\frac{{2k\choose k}^2}{16^k}\right)^2{\bigg /}
n^2{n-1\choose \frac{n-1}{2}}^2.\label{dd2}
\end{align}
\end{lem}
{\noindent \it Proof.}
The above identity possesses, of course, the same automated proof as \eqref{cc5}, and we omit the details.
\qed

{\noindent \it Proof of Theorem \ref{t2}.}
Assume that $p\equiv 3\pmod{4}$. Letting $n=\frac{p-1}{2}$ in \eqref{dd1}, we obtain
\begin{align}
\sum_{k=0}^{\frac{p-1}{2}}\frac{\left(\frac{1}{2}\right)_k\left(\frac{1+p}{2}\right)_k\left(\frac{1-p}{2}\right)_k}{(1)_k^3}=0.
\label{dd3}
\end{align}
Similarly to the proof of \eqref{cc8}, we can show that for $0\le k\le \frac{p-1}{2}$,
\begin{align}
\frac{\left(\frac{1+p}{2}\right)_k\left(\frac{1-p}{2}\right)_k}{(1)_k^2}
\equiv \left(\frac{\left(\frac{1}{2}\right)_k}{(1)_k}\right)^2\left(1-p^2\sum_{j=1}^k\frac{1}{(2j-1)^2}\right)
\pmod{p^4}.\label{dd4}
\end{align}
Substituting \eqref{dd4} into \eqref{dd3} gives
\begin{align}
\sum_{k=0}^{\frac{p-1}{2}}\left(\frac{\left(\frac{1}{2}\right)_k}{k!}\right)^3
\equiv p^2\sum_{k=0}^{\frac{p-1}{2}}\left(\frac{\left(\frac{1}{2}\right)_k}{k!}\right)^3
\sum_{j=1}^k\frac{1}{(2j-1)^2}\pmod{p^4}.\label{dd5}
\end{align}
Furthermore, letting $n=\frac{p-1}{2}$ in \eqref{dd2} and then using \eqref{cc10}, we arrive at
\begin{align}
\sum_{k=0}^{\frac{p-1}{2}}\left(\frac{\left(\frac{1}{2}\right)_k}{k!}\right)^3
\sum_{j=1}^k\frac{1}{(2j-1)^2}\equiv -4^{\frac{p-3}{2}}J_{(p-3)/4}^2{\bigg /}
\left(\frac{p-1}{2}\right)^2{\frac{p-3}{2}\choose \frac{p-3}{4}}^2\pmod{p^2}.\label{dd6}
\end{align}

By \eqref{bb4}, we have
\begin{align}
4^{\frac{p-3}{2}}{\bigg /}
\left(\frac{p-1}{2}\right)^2{\frac{p-3}{2}\choose \frac{p-3}{4}}^2
=\frac{1}{4}\cdot\frac{(1)_{\frac{p-3}{4}}^2}{\left(\frac{1}{2}\right)_{\frac{p+1}{4}}^2}
=\frac{1}{4}\cdot\frac{\Gamma_p\left(\frac{p+1}{4}\right)^2\Gamma_p\left(\frac{1}{2}\right)^2}
{\Gamma_p\left(1\right)^2\Gamma_p\left(\frac{p+3}{4}\right)^2}.\label{dd7}
\end{align}
Since $p\equiv 3\pmod{4}$, by \eqref{bb1} and \eqref{bb2} we obtain
\begin{align*}
&\Gamma_p\left(1\right)^2=1,\\
&\Gamma_p\left(\frac{1}{2}\right)^2=(-1)^{\frac{p+1}{2}}=1,\\
&\Gamma_p\left(\frac{p+3}{4}\right)^2\Gamma_p\left(\frac{1-p}{4}\right)^2=1.
\end{align*}
Substituting the above into \eqref{dd7} yields
\begin{align*}
4^{\frac{p-3}{2}}{\bigg /}
\left(\frac{p-1}{2}\right)^2{\frac{p-3}{2}\choose \frac{p-3}{4}}^2=
\frac{1}{4}\Gamma_p\left(\frac{1+p}{4}\right)^2\Gamma_p\left(\frac{1-p}{4}\right)^2.
\end{align*}
By \eqref{bb5}, we have
\begin{align*}
&\Gamma_p\left(\frac{1+p}{4}\right)\equiv \Gamma_p\left(\frac{1}{4}\right)\left(1+\frac{p}{4}G_1\left(\frac{1}{4}\right)\right)\pmod{p^2},\\[10pt]
&\Gamma_p\left(\frac{1-p}{4}\right)\equiv \Gamma_p\left(\frac{1}{4}\right)\left(1-\frac{p}{4}G_1\left(\frac{1}{4}\right)\right)\pmod{p^2}.
\end{align*}
It follows that
\begin{align*}
\Gamma_p\left(\frac{1+p}{4}\right)\Gamma_p\left(\frac{1-p}{4}\right)
\equiv \Gamma_p\left(\frac{1}{4}\right)^2\pmod{p^2},
\end{align*}
and so
\begin{align}
4^{\frac{p-3}{2}}{\bigg /}
\left(\frac{p-1}{2}\right)^2{\frac{p-3}{2}\choose \frac{p-3}{4}}^2\equiv
\frac{1}{4}\Gamma_p\left(\frac{1}{4}\right)^4\pmod{p^2}.\label{dd8}
\end{align}
Then the proof of \eqref{aa4} follows from \eqref{dd5}, \eqref{dd6} and \eqref{dd8}.
\qed


\begin{thebibliography}{99}
\small \setlength{\itemsep}{-.8mm}

\bibitem{cohen-2007}H. Cohen, Number Theory. Vol. II. Analytic and Modern Tools, Grad. Texts in Math., vol. 240, Springer, New York, 2007.
    
\bibitem{guojw-itsf-2017}V.J.W. Guo, Some generalizations of a supercongruence of van Hamme,
Integral Transforms Spec. Funct. 28 (2017), 888--899.

\bibitem{Long-2011}L. Long, Hypergeometric evaluation identities and supercongruences, Pacific J. Math. 249 (2011), 405--418.

\bibitem{lr-2016}L. Long and R. Ramakrishna, Some supercongruences occurring in truncated hypergeometric series, Adv. Math. 290 (2016), 773--808.

\bibitem{mo-2008}D. McCarthy and R. Osburn, A $p$-adic analogue of a formula of Ramanujan, Arch. Math. (Basel) 91 (2008), 492--504.

\bibitem{mortenson-2008}E. Mortenson, A $p$-adic supercongruence conjecture of van Hamme, Proc. Amer. Math. Soc. 136 (2008), 4321--4328.

\bibitem{os-mc-2009}R. Osburn and C. Schneider, Gaussian hypergeometric series and supercongruences, Math. Comp. 78 (2009), 275--292.

\bibitem{oz-2016}R. Osburn and W. Zudilin, On the (K.2) supercongruence of Van Hamme,
J. Math. Anal. Appl. 433 (2016), 706--711.



\bibitem{schneider-1999}C. Schneider, Symbolic summation assists combinatorics, S\'em. Lothar. Combin. 56 (2007), B56b, 36 pp.

\bibitem{slater-1966}L.J. Slater, Generalized Hypergeometric Functions, Cambridge University Press, Cambridge, 1966.
    
\bibitem{sunzw-2012}Z.-W. Sun, A refinement of a congruence result by van Hamme and Mortenson,
¡¡¡¡Illinois J. Math. 56 (2012), 967--979.

\bibitem{swisher-2015}H. Swisher, On the supercongruence conjectures of van Hamme, Res. Math. Sci. 2 (2015), Art. 18, 21 pp.

\bibitem{vanHamme-1997}L. van Hamme, Some conjectures concerning partial sums of generalized hypergeometric series, $p$-adic functional analysis (Nijmegen, 1996), Lecture Notes in Pure and Appl. Math., vol. 192, Dekker, New York, 1997, 223--236.

\bibitem{whipple-1926}F.J.W. Whipple, On well-posed series, generalised hypergeometric series havinparameters in pairs, each pair with the same sum, Proc. London Math. Soc. 24 (1926), 247--263.

\bibitem{zudilin-2009}W. Zudilin, Ramanujan-type supercongruences, J. Number Theory 129 (2009), 1848--1857.

\end{thebibliography}
\end{document}